# Aspects of Randomization in Infinitely Divisible and Max-Infinitely Divisible Laws.

**S Satheesh**

Telecom Training Centre, BSNL
Trichur – 680 001, **India.**

e-mail - ssatheesh@sancharnet.in

**Abstract.** Here we study certain aspects of randomization in infinitely divisible (ID) and max-infinitely divisible (MID) laws. They generalize ID and MID laws. In particular we study mixtures of ID & MID laws, its relation to random sums & random maximums, corresponding stationary processes & extremal processes and some of their properties. It is shown that they appear as the limit law of ID and MID law when the sample size $N$ is random. We identify a class of probability generating functions for $N$. A method to construct class-L laws is given.

**AMS Classification (2000)**: 60E07, 60F05, 60G50, 60G70, 62E10, 62E20, 62H05, 62N05, 62P05.

**1. Introduction.** Mixtures of infinitely divisible (ID) laws were introduced by Feller (1971, p.573) and studied by Keilson and Steutel (1974). From the angle of summation generalizations of ID laws were considered by: Klebanov, et. al (1984), Sandhya (1991, 1996), Sandhya and Pillai (1999), Mohan, et. al (1993),

Visit: http://ernakulam.sancharnet.in/probstat/

Ramachandran (1997), Gnedenko and Korolev (1996), Klebanov and Rachev (1996), Bunge (1996), Kozubowski and Panorska (1996, 1998), Satheesh (2001a,b), Satheesh, et. al (2002) and Satheesh and Sandhya (2002). In the literature we also have the discussion of max-infinitely divisible (MID) laws parallel to the summation scheme. MID laws were studied by Balkema and Resnick (1977), geometrically MID laws by Rachev and Resnick (1991), G-max-stable laws by Sreehari (1995) and the stability of N-extremes of continuous and discrete laws by Satheesh and Nair (2002b).

The investigation presented here has been motivated by (i) the $\varphi$-ID laws of Satheesh (2001b) (ii) the G-max-stable laws Sreehari (1995) (iii) mixtures of Poisson laws studied in Satheesh and Nair (2002a) (iv) Feller's proof of Bernstein's Theorem and (v) Gnedenko's (1982) transfer theorems for random sums and maximums.

**2. Results.** Treating a parameter in a probability model is usually termed as randomization and such schemes generate new models. Schemes of interest here are (i) mixtures, (ii) randomized operational time in stochastic processes and (iii) random sums. We will discuss certain implications of these three schemes on ID and MID laws.

From Feller (1971, p.557) we know that a characteristic function (CF) $\omega$ is ID iff $\omega^s$ is a CF for every $s>0$ and further $\omega$ does not have a real zero. Similarly a d.f $H$ is MID iff $H^s$ is a d.f for every $s>0$. Since all d.f $H$ on $\boldsymbol{R}$ are MID discussion of MID laws and their generalizations are relevant only for d.fs on $\boldsymbol{R}^d$ for $d \geq 2$ integer. Further the operations are to be taken component wise. Thus in the context of MID laws all d.fs are defined on $\boldsymbol{R}^d$ for $d \geq 2$ integer and the set $\{x \in \boldsymbol{R}^d : H(x)>0\}$ will be denoted by $\{H>0\}$. From Balkema and Resnick (1977) we have: $F$ is MID iff $\{H>0\}$ is a rectangle. Obviously the discussion is on $\{F>0\}$. Let us now consider randomizing the parameter $s>0$.

Let $Z>0$ be a r.v with d.f $G$ and Laplace transform (LT) $\varphi$.

**Theorem 2.1a** If the parameter $s$ in the CF $\omega^s$ of an ID law is randomized assigning the d.f $G$, the resultant distribution will have the CF $\varphi\{-\log \omega\}$.

*Proof.* The conclusion follows from the relation: $\int_0^\infty \exp\{s \log \omega\} dG(s) = \varphi\{-\log \omega\}$.

These models with CF $\varphi\{-\log \omega\}$ will be called $\varphi$-mixtures of ID laws. Similarly, we can describe $\varphi$-mixtures of MID laws as follows.

**Theorem 2.1b** If the parameter $s$ in the d.f $H^s$ of a MID law is randomized assigning the d.f $G$, the resultant distribution will have the d.f $\varphi\{-\log H\}$.

Randomizing the operational time of a stochastic process has been discussed in Feller (1971, p.345). In $\varphi$-mixtures of ID and MID laws this means:

**Theorem 2.2a** For every $\varphi$-mixture of ID law with CF $h = \varphi\{-\log \omega\}$, the CF $h^t$ represents a stochastic process $\{X(T(t)), t>0\}$ given $X(0) = 0$; that is subordinated to the process $\{X(s), s>0\}$ having stationary and independent increments with CF $\omega^s$, by the directing process (operational time) $\{T(t), t>0\}$ with LT $\varphi^t$.

*Proof.* Follows from Feller (1971, p.345 & p.573) and Theorem 2.1a.

**Theorem 2.2b** Corresponding to every d.f $F$ that is a $\varphi$-mixture of a MID law with d.f $H$ there exists a multivariate extremal process $\{Y(t), t>0\}$ governed by this MID law and an independent r.v $Z$ with d.f $G$ and LT $\varphi$ such that $F(x) = P\{Y(Z) \leq x\}$.

*Proof.* Let $\{Y(t), t>0\}$ be the extremal process governed by the MID law with d.f $H$. That is; $P\{Y(t) \leq x\} = \exp\{t \log H\}$.

Hence; $P\{Y(Z) \leq x\} = \int_o^\infty \exp\{t \log H\} dG(t) = \varphi\{-\log H\} = F(x)$.

Sandhya and Satheesh (1996), Ramachandran (1997) and Satheesh, et. al (2002) have discussed the relation between random sum stability and class-L ($L$) laws. Here we discuss this from the angle of a φ-mixture and generalize Theorem.2.1 (part (iii), (iv) and (v)) of Ramachandran (1997) and Theorems 2.3, 2.4 and 2.5 regarding generalized Linnik laws of Satheesh, et. al (2002).

**Property 2.1** A φ-mixture of a strictly stable law is in $L$ if $\varphi \in L$.

*Proof.* If $\varphi \in L$ then there is another LT $\varphi_c$ such that

$$\varphi(s) = \varphi(cs) \, \varphi_c(s), \; s>0 \text{ for every } 0<c<1.$$

Now, considering the φ-mixture of a strictly stable law we have

$$\varphi(\psi(t)) = \varphi(c\psi(t)) \, \varphi_c(\psi(t)),$$

where $\psi(t) = \lambda|t|^\alpha e^{-i\theta \text{sgn}(t)}$, $\lambda>0$, $|\theta| \leq \min(\pi\alpha/2, \pi-\pi\alpha/2)$, $0<\alpha \leq 2$.

Showing that the φ-mixture of a strictly stable law is in $L$.

**Theorem 2.3** The generalized Linnik $(\lambda, \alpha, \theta, \nu)$ law with CF $\{1+ \lambda|t|^\alpha e^{-i\theta \text{sgn}(t)}\}^{-\nu}$, $\lambda>0$, $0<\alpha \leq 2$, $|\theta| \leq \min(\pi\alpha/2, \pi-\pi\alpha/2)$, $\nu>0$, is in $L$.

*Proof.* This CF can be obtained as the gamma$(\lambda,\nu)$-mixture of strictly stable law. Since gamma$(\lambda,\nu)$ law is in $L$ this mixture is also in $L$.

**Remark 3.1** An implication of this is that this class of distributions is absolutely continuous and unimodal. Hence the generalized Linnik$(\alpha,\theta,\nu)$ law (the above CF with $\lambda = 1$) of Erdogan and Ostrowski (1998) is in $L$ and they are absolutely continuous and unimodal. This is a much simpler line of argument than the one given by them.

**Remark 3.2** Having membership in $L$ or not is important from a modeling perspective as such distributions can model ARMA processes, Gaver and Lewis (1980). Property.2.1 also enables one to construct $L$ laws.

**Note.** Property 2.1 suggests the notion of $L$-analogue for maximums and that of Theorem 2.3.

The following two lemmas are ramifications of Feller's proof (Feller, 1971, p.440) of Bernstein's theorem and were first proved in Satheesh (2001b). Together they provide a class of discrete laws that can be used in the transfer theorems for sums and extremes.

**Lemma 2.1** $\wp_\varphi(s) = \{P_\theta(s) = s^j \varphi\{(1-s^k)/\theta\}, 0<s<1, j \geq 0 \text{ \& } k \geq 1 \text{ integer and } \theta > 0\}$ describes a class of probability generating functions (PGF) for any given LT $\varphi$.

*Proof.* From Feller (1971, p.440) it follows that $\varphi\{(1-s)/\theta\}$ is a PGF for any given LT $\varphi$. Let $N$ be the r.v with PGF $\varphi\{(1-s)/\theta\}$ and the r.v $X$ is such that $P(X=k) = 1$, $k \geq 1$ integer and independent of $N$. If $X_1, X_2, \ldots$ are independent and identically distributed as $X$, then $S_N = X_1 + \ldots + X_N$ has the PGF $\varphi\{(1-s^k)/\theta\}$. Now the function $s^j \varphi\{(1-s^k)/\theta\}$, $j \geq 0$ corresponds to the PGF of a distribution with the support shifted $j$ integers to the right. (This proof brings out the structure of the distribution. For a short analytic proof see Satheesh (2001b)).

**Lemma 2.2** Given a r.v $U$ with LT $\varphi$, the integer valued r.vs $N_\theta$ with PGF $P_\theta$ in the class $\wp_\varphi(s)$ described in Lemma.2.1 satisfy $\theta N_\theta \xrightarrow{d} kU$ as $\theta \to 0$.

*Proof.* The LT of $\theta N_\theta$ is $e^{-vj\theta} \varphi\{(1-e^{-vk\theta})/\theta\}$, $v>0$ and by Feller (1971, p.440),

$$\underset{\theta \to 0}{Lt} \{e^{-vj\theta} \varphi((1-e^{-vk\theta})/\theta)\} = \varphi(kv).$$

Now, let $\{X_{\theta,j}; \theta \in \Theta, j \geq 1\}$ be i.i.d r.vs with CF $g_\theta$ and for $k \geq 1$ integer, set $S_{\theta,k} = X_{\theta,1} + \ldots + X_{\theta,k}$. Similarly consider i.i.d random vectors $\{\mathbf{Y}_{\theta,j}; \theta \in \Theta, j \geq 1\}$ in $\mathbf{R}^d$ for $d \geq 2$ integer with d.f $G_\theta$ and $\mathbf{M}_{\theta,k} = \vee \{\mathbf{Y}_{\theta,1}, \ldots, \mathbf{Y}_{\theta,k}\}$. In the sequel we will consider a $\{\theta_n \in \Theta\}$ and $\{k_n; n \geq 1$ of positive integers$\}$ such that as $n \to \infty$, $\theta \downarrow 0$ through $\{\theta_n\}$ and $S_{\theta,k_n} \xrightarrow{d} U$ so that $U$ is ID and $\mathbf{M}_{\theta,k_n} \xrightarrow{d} \mathbf{V}$ so that $\mathbf{V}$ is MID. Also in an $N_\theta$-sum of $\{X_{\theta,j}\}$ ($N_\theta$-max of $\{\mathbf{Y}_{\theta,j}\}$), $N_\theta$ and $X_{\theta,j}$ ($\mathbf{Y}_{\theta,j}$) are assumed to be independent for each $\theta \in \Theta$. The next two results follow from Lemma 2.1, 2.2 and the transfer Theorems for sums and maxs.

**Theorem 2.4a** The limit law of $N_\theta$-sums of i.i.d r.vs $\{X_{\theta,j}\}$ as $n \to \infty$, where the PGF of $N_\theta$ is a member of $\wp_\varphi(s)$, is necessarily a $\varphi$-mixture of an ID law.

**Theorem 2.4b** The limit law of $N_\theta$-maxs of i.i.d random vectors $\{\mathbf{Y}_{\theta,j}\}$ as $n \to \infty$, where the PGF of $N_\theta$ belongs to $\wp_\varphi(s)$, is necessarily a $\varphi$-mixture of a MID law.

From these two theorems one can conceive the notions of attraction and partial attraction for $N_\theta$-sums and $N_\theta$-maxs when the PGF of $N_\theta$ belongs to $\wp_\varphi(s)$. Here we will denote the PGFs by $P_n$ -corresponding to $\theta_n = 1/n$, $\{n\}$ the sequence of positive integers and $P_{n_m}$ -corresponding to $\theta_m = 1/n_m$, $\{n_m\}$ a subsequence of $\{n\}$.

**Definition 2.1a** A CF $g(t)$ belongs to the domain of $\varphi$-attraction (D$\varphi$-A) of the CF $f(t)$ if there exist sequences of real constants $a_n = a(\theta_n) > 0$ and $b_n = b(\theta_n)$ such that with $g_n(t) = g(t/a_n) \exp(-itb_n)$; $\underset{n \to \infty}{Lt} P_n\{g_n(t)\} = f(t)$, $\forall t \in \mathbf{R}$, and $g(t)$ belongs to the domain of partial $\varphi$-attraction (DP$\varphi$-A) of $f(t)$ if; $\underset{m \to \infty}{Lt} P_{n_m}\{g_m(t)\} = f(t)$, $\forall t \in \mathbf{R}$. Similarly:

**Definition 2.2b** A d.f $G$ in $\mathbf{R}^d$ belongs to the D$\varphi$-MA ($\varphi$-max-attraction) of the d.f $F$ with non-degenerate marginal distributions, if there exists normalizing constants $a_{i,n} = a_i(\theta_n) > 0$ and $b_{i,n} = b_i(\theta_n)$, $1 \leq i \leq d$, such that with $G_n(\mathbf{y}) = G(a_{i,n} y_i + b_{i,n}, 1 \leq i \leq d)$; $\underset{n \to \infty}{Lt}\ P_n\{G_n\} = F$ and $G$ belongs to the DP$\varphi$-MA of $F$ if; $\underset{m \to \infty}{Lt}\ P_{n_m}\{G_m\} = F$.

Setting $\theta_m = 1/n_m$ ($\theta_n = 1/n$) and invoking the transfer theorems (or making these prescriptions in Theorems 2.4a and 2.4b) we have: if the CF $g(t)$ is in the DPA (DA) of an ID (stable) law, then it is in the DP$\varphi$-A (D$\varphi$-A) of the corresponding $\varphi$-mixture of ID (stable) law. Clearly the analogous results hold good for $N_\theta$-maxs and $\varphi$-mixture of MID (max-stable) laws. Satheesh (2001b) and Satheesh and Sandhya (2002) have investigated and shown that the converses are also true.

**Concluding Remarks.** Generalizations of ID and geometrically-ID laws were considered by: Sandhya (1991, 1996), discussing two examples of non-geometric laws for $N$ but the description was not constructive. The description of N-ID laws (with CF $\varphi\{\psi\}$, where $\exp\{-\psi\}$ is the CF of an ID law) by Gnedenko and Korolev (1996), Klebanov and Rachev (1996) & Bunge (1996), are based on the assumption that the PGF $\{P_\theta,\ \theta \in \Theta\}$ of $N$ formed a commutative semi-group. Satheesh (2001b) showed that this assumption is not natural since it captures the notion of stable rather than ID laws. Also, it rules out any $P_\theta$ having an atom at the origin. A glaring situation is that this theory cannot approximate negative binomial sums by taking $\varphi$ as gamma as one expects. In the case of exact distributions Satheesh, et. al (2002) showed that when $\varphi$ is gamma, $N_\theta$ is Harris. Kozubowski and Panorska (1996, 1998) discussed $\nu$-stable laws with CF $\varphi\{\psi\}$ ($\exp\{-\psi\}$ being a stable CF) approximating $N_\theta$-sums, assuming that $N_\theta \xrightarrow{p} \infty$ as $\theta \downarrow 0$ and it could handle negative binomial sums. But identifying those $N_\theta \xrightarrow{p} \infty$ has not been done.

Satheesh (2001b), motivated by these observations and to overcome them introduced φ-ID laws with CF $\varphi\{\psi\}$ approximating $N_\theta$-sums where the PGF of $N_\theta$ is a member of $\mathscr{P}_\varphi(s)$ (that of Lemma 2.1) without requiring $N_\theta \xrightarrow{p} \infty$. For this class of PGFs $\mathscr{P}_\varphi(s)$, he proved a result, that is analogous to the Theorem 4.6.5 in Gnedenko and Korolev (1996, p.149), and discussed φ-attraction and partial φ-attraction. His results generalize those on attraction and partial attraction in geometric sums in Sandhya (1991), Sandhya and Pillai (1999), Mohan et al (1993), Ramachandran (1997), and those for N-sums in Gnedenko and Korolev (1996). Satheesh and Sandhya (2002) shows that when the PGF of $N_\theta$ is a member of $\mathscr{P}_\varphi(s)$ a N&S condition for the convergence of $N_\theta$-sums of $\{X_{\theta,j}\}$ to a φ-ID law is $(1 - g_\theta(t))/\theta \to \psi(t)$ as $\theta \downarrow 0$ through a $\{\theta_n\}$, where $\psi(t)$ is a continuous function. This result generalizes Theorem.1.1 of Feller (1971, p.555). Analogous results for random vectors and MID laws are also given therein, describing the notion of partial φ-max-attraction for φ-MID laws with d.f $\varphi\{-\log H\}$ (of Theorem.2.1b) thus generalizing the results for the geometric case given in Rachev and Resnick (1991).

**Acknowledgement.** The author is thankful to the referee for the suggestions.